\documentclass[12 pt]{article}
\usepackage{amssymb,amsmath,amsfonts,amsthm}
\renewcommand{\P}{\ensuremath{\mathbb {P}}}

\newcommand\R{{\mathbb {R}}}

\newcommand\E{{\mathbb {E}}}

\newcommand\sigmas{\sigma_{\sup}}

\newtheorem {Theorem}{Theorem}[section]

\newtheorem {Proposition}[Theorem]{Proposition}

\theoremstyle{definition}
\newtheorem{Definition}{Definition}[section]

\newtheorem{Remark}{Remark}[section]

\newcommand\F{{\mathcal {F}}}

\newcommand\beq{\begin{equation}}
\newcommand\eeq{\end{equation}}

\begin{document}
\title{Large and moderate deviations for the left random 
walk on $GL_d({\mathbb R})$}

\date{}

\author{Christophe Cuny$^a$, J\'er\^ome Dedecker$^b$ and Florence Merlev\`ede$^c$}

\maketitle

{\abstract{Using martingale methods, we obtain some upper bounds for large and moderate deviations of products of independent and identically distributed elements of
$GL_d({\mathbb R})$.  We investigate all the possible moment conditions, from super-exponential moments  to weak moments of 
order $p>1$,  to get a complete picture of the situation. 
We also prove a moderate deviation principle under an 
appropriate tail condition.} }

\bigskip

\small{

$^a$ Universit\'e de la Nouvelle-Cal\'edonie, Equipe ERIM, et CentraleSupelec, Laboratoire MAS.

Email: cuny@univ-nc.nc \bigskip

$^b$ Universit\'e Paris Descartes, Sorbonne Paris Cit\'e,  Laboratoire MAP5 (UMR 8145).

Email: jerome.dedecker@parisdescartes.fr \bigskip

$^c$ Universit\'{e} Paris-Est, LAMA (UMR 8050), UPEM, CNRS, UPEC.

Email: florence.merlevede@u-pem.fr}

\section{Introduction}
Let $(\Omega,\F,\P)$ be a probability space and $(Y_n)_{n\ge 1}$  be independent and identically distributed random variables on $(\Omega,\F,\P)$ taking values 
in $G:=GL_d(\R)$, $d\ge 2$ (the group of invertible $d$-dimensional real matrices),
with common distribution $\mu$.
Denote by $\Gamma_\mu$ 
the closed semi-group generated by the support of $\mu$.
Let $\|\cdot\|$ be the euclidean norm on $\R^d$, and 
for every
$g\in GL_d(\R)$, let $\|g\|:=\sup_{\|x\|=1}\|gx\|$.

In all the paper, we assume that 
$\mu$ is \emph{strongly irreducible}, i.e. that no proper finite union of subspaces of $\R^d$ are invariant by $\Gamma_\mu$ and that 
it is \emph{proximal}, i.e. that there exists a matrix in $\Gamma_\mu$ 
admitting a unique (with multiplicity one) eigenvalue with maximum modulus. 

For such a measure $\mu$, 
it is known  that there exists a unique invariant measure 
$\nu$ on the projective space $X:=P_{d-1}(\R)$ 
 (see for instance Theorem 3.1 of \cite{BL}) in the following sense:  for any bounded Borel function $h$ from $X$ to ${\mathbb R}$
\begin{equation} \label{inv22-07}
  \int_X h(u) \nu(du) = \int_G \int_X  h(g \cdot  u) \mu(dg) \nu(du) \, .
\end{equation}
Moreover, if 
\begin{equation}\label{moment-1}
\int_G \log N(g) \mu(dg)<\infty \, , \quad 
\text{where $N(g):=\max(\|g\|,\|g^{-1}\|)$},
\end{equation}
then  (see for instance Corollary 3.4 page 54 of \cite{BL} or 
Theorem 3.28 of \cite{BQ-book}),
for every $x \in S^{d-1}$, 
\begin{equation}\label{stronglaw}
\frac{\log \| Y_n \cdots Y_1 x \|}{n}\underset{n\to +\infty}\longrightarrow \lambda_\mu= \int_G\int_X \sigma(g,u) \mu(dg)\nu(du)\quad \text{almost surely}, 
\end{equation}
where 
$$
\sigma(g,\bar x)= 
\log \left (\frac{\|g\cdot x\|}{\|x\|} \right )\, \quad 
\text{for 
$g\in GL_d(\R)$ and  $x\in \R^d-\{0\}$, $\bar x$ denoting the class 
of $x$ in $X$}  \, .
$$
Note  that  the function $\sigma$ defined above is a 
$cocycle$, in the following sense:
\begin{equation}\label{cocycle-prop}
\sigma(gg',u)=\sigma(g,g'\cdot u)+\sigma(g',u)\, \quad 
\text{ for any $g,g'\in G$ and 
$u\in X$}\, .
\end{equation}

Let $A_k=Y_k \cdots Y_1$ for $k\geq 1$ and 
$A_0=\mathrm{Id}$. 
In this paper
we wish to study the asymptotic behavior of 
\begin{equation}\label{modev}
\sup_{\|x\|=1} {\mathbb P} \left ( 
\max_{1 \leq k \leq n} \left |\log \| A_k x \| - k \lambda_\mu \right |
>n^\alpha y \right ) \, ,
\end{equation}
when $\alpha \in (1/2, 1]$, under stronger moment conditions
on $\log N(Y_1)$ than \eqref{moment-1}. This is a way to study rates
of convergence in the strong law \eqref{stronglaw}. 
In the probabilistic terminology, the case $\alpha \in (1/2,1)$
corresponds to the \emph{moderate deviation} regime, and 
$\alpha=1$ to the \emph{large deviation} regime. 

The case $\alpha=1/2$ corresponds to the normalization of  the central limit theorem.
In that case, the asymptotic behavior of \eqref{modev} is due 
to Benoist and Quint \cite{BQ} as soon as $\log N(Y_1)$ has a moment of order 2
(note that Benoist and Quint do not deal with the maximum in 
\eqref{modev}, but their method also applies in that case, 
see also Theorem 1(ii) in \cite{CDJ}). A previous result is due to Jan \cite{Jan} under a moment of order $2+\epsilon$, 
$\epsilon >0$. 

In this paper, we shall give precise informations on the rate of convergence to 0
(as $n \rightarrow \infty$) of \eqref{modev} when 
$\alpha \in (1/2,1]$, under various moment conditions on the random variable $\log N(Y_1)$: sub or super-exponential moments in Section \ref{subsuper}, 
weak moments of order
$p>1$ in Section \ref{weakMom}, and strong moments of order $p\geq 1$
in Section \ref{strongMom}. 

In Section \ref{sec:MDP} we shall give a moderate
deviation principle for the process 
$$\left \{ n^{-\alpha}(\log \| A_{[nt]} x \|-[nt] \lambda_{\mu}), t \in [0,1]  \right \}$$  when $\log N(Y_1)$ satisfies  Arcones's tail
 condition \cite{Arc} (which is true under an 
appropriate sub-exponential moment condition). In Section
\ref{strongMom} we obtain some results in the spirit of 
Baum and Katz \cite{BK} which complement the results
on complete convergence obtained in \cite{BQ} in the case 
$\alpha=1$. When  $\log N(Y_1)$ has a strong moment of
order $p \in (1,2)$ and $\alpha=1/p$, this gives the
rate $n^{(p-1)/p}$ in the strong law of large numbers, 
which was proved by another method in \cite{CDJ}, Theorem 1(i). 

\smallskip

All along the paper, the following notations will be used: let ${\mathcal F}_0 = \{ \emptyset, \Omega \}$ and  ${\mathcal F}_k = \sigma (Y_1, \ldots, Y_k)$,  for any $k \geq 1$. For any $x \in S^{d-1}$, define $X_{0,x}=x$  and $X_{n,x}=\sigma(Y_n, A_{n-1} x)$  for  $n \geq 1$. 
With these notations,   for any $x \in S^{d-1}$ and any positive integer $k$,
\begin{equation} \label{trivialF}
\log \| A_k x \|=  \sum_{i=1}^k X_{i,x} \, .
\end{equation}
The equality \eqref{trivialF} follows easily from the fact 
that $\sigma$ is a cocycle ($i.e.$ \eqref{cocycle-prop} holds). 
In Section \ref{Sec:GC} we shall present some extensions of our results
to general cocycles, in the spirit of Benoist and Quint \cite{BQ}.

\section{The case of (sub/super) exponential moments}
\label{subsuper}
\setcounter{equation}{0}
\subsection{Upper bounds for large deviations}
Let $r>0$. In this subsection, we assume that 
\begin{equation}\label{exp}
  \int_G e^{\delta (\log N(g))^r} \mu(dg) < \infty\, , \quad 
  \text{for some $\delta >0$.}
\end{equation}

We first consider the case $r \geq 1$.  In that case, using the spectral gap property, Le Page \cite{L} proved the following large deviation principle: there exists 
a positive constant $A$ such that, for any $y \in (0,A)$, 
\begin{equation}\label{LD0}
\lim_{n \rightarrow \infty} \frac 1 n 
 \log 
 {\mathbb P} \left ( 
 \left |\log \| A_n x \| - n \lambda_\mu \right |
>n y \right ) =\phi(y) \, .
\end{equation}
Of course, this is the best possible result for $y \in (0,A)$. However, it does not give any information for large values of $y$, and the rate 
function $\phi$ is not explicit (in particular, one cannot easily describe the behavior of $\phi$ when $r$ varies in $[1, \infty)$). 

The following result, which is obtained $via$ a completely different method, can be seen as a 
complementary result of \eqref{LD0}. It gives an explicit (up to a constant) upper bound for $\phi$ when $y \in (0,A)$, and this upper bound is valid for any $y>0$. 
In particular, we can see the qualitative change in the behavior of large deviations for large $y$ when $r$ varies in $[1, \infty)$. 
\begin{Theorem}\label{LDsuperexp}
Assume that \eqref{exp} holds for some $r \geq 1$. 
Then there exists a positive constant $C$ such that, for any 
$y>0$,
\begin{equation}\label{LD1}
\limsup_{n \rightarrow \infty} \frac 1 n 
 \log 
\sup_{\|x\|=1} {\mathbb P} \left ( 
\max_{1 \leq k \leq n} \left |\log \| A_k x \| - k \lambda_\mu \right |
>n y \right ) \leq
-C \left( y^2 {\bf 1}_{y \in (0,1)}+ y^r{\bf 1}_{y \geq 1}\right )\, .
\end{equation}
\end{Theorem}

 For $r \in (0,1)$, there is no such result as \eqref{LD0}. Instead, one can prove:
\begin{Theorem}\label{LDsousexp}
 Assume that \eqref{exp} holds for some  $r \in (0,1)$.
Then there exists a positive constant $C$ such that, 
for  any $y>0$, 
\begin{equation}\label{LD2}
\limsup_{n \rightarrow \infty} \frac {1} {n^r} \log 
\sup_{\|x\|=1}
{\mathbb P} \left ( \max_{1 \leq k \leq n}
\left | \log \| A_k x \|  -k \lambda_\mu \right |>n y 
\right ) \leq
-C y^r \, .
\end{equation}
\end{Theorem}

\noindent{\bf Proofs of Theorems \ref{LDsuperexp} and \ref{LDsousexp}.} 
Since $\int (\log N(g))^2 \mu (dg) <\infty $, we infer  from the equality (3.9) in \cite{BQ} that, for any $x \in S^{d-1}$,
\begin{equation}\label{MartCob1}
X_{k,x}- \lambda_\mu 
= D_{k,x} + \psi (A_{k-1}x)- \psi (A_k x) \, ,
\end{equation}
where $\psi$ is a bounded function and $D_{k,x}$ is 
${\mathcal F}_k$-measurable and such that 
${\mathbb E}(D_{k,x} | {\mathcal F}_{k-1})=0$. The decomposition \eqref{MartCob1} 
is  called a \emph{martingale-coboundary} decomposition. Such a decomposition 
has been used for the first time in the paper \cite{Gordin} by Gordin (see also \cite{GL}).

\medskip

Starting from \eqref{trivialF} and \eqref{MartCob1}, for any $x \in S^{d-1}$ and any positive integer $k$,
\begin{equation}\label{MartCob2}
\log \| A_k x \|-k \lambda_\mu = M_{k,x} + \psi(x)- 
\psi (A_{k} x) \, ,
\end{equation}
where $M_k(x)= D_{1,x} + \cdots +D_{k,x}$ is a martingale 
adapted to the filtration ${\mathcal F}_k$. Clearly, 
since $|\psi(x)- 
\psi (A_{k} x)| \leq 2 \| \psi \|_\infty$, it is equivalent to 
prove \eqref{LD1} and \eqref{LD2} for $M_{k,x}$ instead of 
$(\log \| A_k x \|-k \lambda_\mu )$. 

To do this, we first note that if \eqref{exp} holds, then 
\begin{equation}\label{boundXexp}
\left \| {\mathbb E} \left ( e^{ \delta |X_{k,x}|^r} \Big | 
\mathcal{F}_{k-1} \right ) \right \|_\infty 
=
\left \| \int  e^{ \delta |\sigma(g, A_{k-1}x)|^r}  \mu (dg) 
\right \|_\infty 
\leq 
\int_G e^{\delta (\log N(g))^r} \mu(dg) < \infty \, .
\end{equation}
Using again that $\psi$ is bounded  we infer from \eqref{MartCob1} and \eqref{boundXexp}  that there
 exists
a constant $K$ such that
\begin{equation}\label{boundM}
\sup_{\|x\|=1} 
\left \| {\mathbb E} \left ( e^{ \delta |D_{k,x}|^r} \Big | 
\mathcal{F}_{k-1} \right ) \right \|_\infty < K \, ,
\end{equation}
for any positive integer $k$. 

Starting from \eqref{boundM}, it remains to apply known results
to the martingale $M_k(x)$. 

To prove \eqref{LD1} (case $r\geq 1$), we apply Theorem 1.1 of  \cite{LW}, which implies that there exists a positive constant $c$
such that, for any $y>0$,
\begin{equation}\label{LW}
\sup_{\|x\|=1} {\mathbb P} \left ( \max_{1 \leq k \leq n} |M_{k,x}|> n y \right )
 \leq 2 \exp \left ( - n c \left ( y^2 {\bf 1}_{y \in (0,1)}+ 
 y^r{\bf 1}_{y \geq 1} \right ) \right ) \, .
\end{equation}
The upper bound \eqref{LD1} follows directly from \eqref{LW}. 
Note that a direct application of  Theorem 1.1 of \cite{LW} gives 
us \eqref{LW} without the maximum over $k$. However, a careful
reading of the proof reveals that one can take the maximum
over $k$. The only argument that should be added to the proof  is Doob's maximal inequality for non-negative submartingales,
which implies that   
$$
{\mathbb E}\left ( e^{\lambda \max_{1 \leq k \leq n} M_{k,x}}
\right ) \leq {\mathbb E}\left (e^{ \lambda  M_{n,x}}
\right ) ,  \ \text{for any $\lambda >0$.}
$$

To prove \eqref{LD2} (case $r \in (0,1)$), we apply Theorem 2.1 of \cite{Fan} (see
also the proof of Proposition 3.5 of \cite{DF}) and more precisely 
the upper bound (13) in \cite{Fan}, which implies that there exist
two positive constants $c_1$ and $c_2$ such that 
\begin{equation}\label{F}
\sup_{\|x\|=1} {\mathbb P} \left ( \max_{1 \leq k \leq n} |M_{k,x}|> n y \right ) \leq 4 \exp \left ( - c_1 (ny)^r  \right ), \
\text{for any $y>c_2 n^{-(1-r)/(2-r)}$.}
\end{equation}
The upper bound \eqref{LD2} follows directly from \eqref{F}. $\Diamond$

\subsection{A moderate deviation principle}\label{sec:MDP}
Let $(b_n)_{n\geq 0}$ be a sequence of positive numbers satisfying the following regularity conditions: 
\begin{equation} \label{condseq}
\text{The functions $f(n)=\frac{n^2}{b_n^2}$ and $g(n)=b_n^2$ are strictly increasing to infinity, and $\displaystyle \lim_{n \rightarrow \infty} \frac{n}{b_n^2} = 0$.}
\end{equation}

For $x \in S^{d-1}$, let
$$
Z_{n,x}=
\left \{ \frac{\log \| A_{[nt]} x \|-[nt] \lambda_{\mu}}{b_n}, t \in [0,1] \right \} \, .
$$
The process $Z_{n,x}$ takes values in the space $D([0,1])$ equipped with the usual Skorokhod topology. 
The following functional moderate deviation principle holds:

\begin{Theorem}\label{MDPth}
Let $(b_n)_{n\geq 0}$ be a sequence of positive numbers  satisfying the condition \eqref{condseq}. Assume  $\int (\log N(g))^2 \mu (dg) <\infty $ and 
\begin{equation}\label{ARC}
 \limsup_{n \rightarrow \infty} 
  \frac{n}{b_n^2}\log n \mu \left \{ \log N > b_n \right \} 
  = - \infty \, .
\end{equation}
Then, for any $x \in S^{d-1}$, 
$
n^{-1} {\mathbb E}((\log \|A_n x \|- n \lambda_\mu)^2) 
\rightarrow V
$ as $n\rightarrow \infty$, where $V$ does not depend  on $x$. 
Moreover, for any 
Borel set $\Gamma \subset D([0,1])$,
\begin{equation}\label{fmdp}
- \inf_{ \varphi  \in \Gamma^\circ}  I_V(\varphi) \leq \liminf_{n \rightarrow \infty} 
\frac{n}{b_n^2} \log \inf_{\|x\|=1}
{\mathbb P} \left (
 Z_{n,x} \in \Gamma 
\right ) 
\leq 
\limsup_{n \rightarrow \infty} 
\frac{n}{b_n^2}  \log \sup_{\|x\|=1}
{\mathbb P} \left ( Z_{n,x}  \in \Gamma 
\right )  \leq - \inf_{ \varphi  \in \bar \Gamma}  I_V(\varphi) \, ,
\end{equation}
where
$$
I_V (h) = \frac{1}{2 V} \int_0^1 \left (
h'(u) \right )^2 du 
$$
if simultaneously $V>0$, $h(0)=0$ and $h$ is absolutely 
continuous, and $I_V (h)= + \infty$ otherwise.
\end{Theorem}

\begin{Remark} Let $(b_n)_{n\geq 0}$ be a sequence of positive numbers  satisfying \eqref{condseq}. If $(X_i)_{i \geq 1}$ is a sequence of independent and identically 
distributed random variables, Arcones \cite{Arc} proved 
that the functional  moderate deviation principle holds provided $\E (X_1^2) < \infty$ and 
\begin{equation}\label{ARC2}
 \limsup_{n \rightarrow \infty} 
  \frac{n}{b_n^2}\log n {\mathbb P}(|X_1|>b_n)
  = - \infty \, .
\end{equation}
Moreover, he showed that  condition \eqref{ARC2} is also necessary for the moderate deviation principle. 
Note that our condition \eqref{ARC} is exactly Arcones's tail condition for the random variable $\log N(Y_1) $. 
When $b_n= n^\alpha$ with $\alpha \in (1/2,1)$, the tail condition  \eqref{ARC}
is true if
$$
\mu \left \{ \log N > x \right \} \leq  e^{- x^\beta a(x) }  \, ,
$$
for  $\beta= 2-(1/\alpha)$ and a function $a$ such that
 $a(x) \rightarrow \infty $ as $x \rightarrow \infty$ (note that $\beta \in (0,1)$, so only a sub-exponential moment is needed for $\log N(Y_1)$).
\end{Remark}

\begin{Remark} Applying the contraction principle, 
Theorem \ref{MDPth} implies in particular that, for any 
Borel set $\Gamma \subset {\mathbb R}^+$,
\begin{multline*}
- \inf \left \{ \frac{y^2}{2V}, y \in \Gamma^\circ \right \} \leq \liminf_{n \rightarrow \infty} 
\frac{n}{b_n^2} \log \inf_{\|x\|=1}
{\mathbb P} \left ( 
 \frac { \max_{1 \leq k \leq n}
\left | \log \| A_k x \|  -k \lambda_\mu \right |}
{ b_n} \in \Gamma 
\right ) \\
\leq 
\limsup_{n \rightarrow \infty} 
\frac{n}{b_n^2} \log \sup_{\|x\|=1}
{\mathbb P} \left (
 \frac { \max_{1 \leq k \leq n}
\left | \log \| A_k x \|  -k \lambda_\mu \right |}
{ b_n} \in \Gamma 
\right )  \leq - \inf \left \{ \frac{y^2}{2V}, y \in \bar \Gamma \right \}\, .
\end{multline*}
Note that a partial result in this direction has been obtained in \cite{BQ-book}, Proposition 11.12. In that Proposition, the authors proved a moderate deviation principle for 
$(\log \| A_n x \|  -n \lambda_\mu)$ and the collection of
open intervals, 
under an exponential moment for $\log N(Y_1)$. However, their result is stated in a more general framework than ours (see Section \ref{Sec:GC} of the present paper  for an extension of Theorem  \ref{MDPth} to general cocycles). 
\end{Remark}

\noindent{\bf Proof of Theorem \ref{MDPth}.}
Since $\int (\log N(g))^2 \mu (dg) <\infty $, the decomposition \eqref{MartCob2} holds, and 
it is equivalent to 
prove \eqref{fmdp}  for the process
$$
\tilde Z_{n,x}=
\left \{ \frac{M_{[nt],x}}{b_n}, t \in [0,1] \right \} 
$$
 instead of 
$Z_{n,x}$.  Now, by a standard argument, to get the  result uniformly with respect to $x \in S^{d-1}$ in \eqref{fmdp}, it suffices to prove the functional moderate deviation principle  for the process
$\tilde Z_{n,x_n}$, where $(x_n)_{n \geq 1}$ is any sequence of points in $S^{d-1}$. 

The result will follow from the next proposition, which is a triangular version of Theorem 1 in \cite{Dj}. This proposition  is in fact a corollary of a more general result for triangular arrays of 
martingale differences which can be deduced from Puhalskii's results and Worms's paper (see \cite{Pu} and \cite{Wo}). We refer to Theorem \ref{PWDFMDP} of  the Appendix 
for a complete statement and some elements of proof. 

Before giving the statement of  this proposition, we need more notations. Assuming \eqref{condseq}, we can  construct the strictly increasing
continuous function $f(x)$ that is formed by the line segments from $%
(n,f(n)) $ to $(n+1,f(n+1))$. Similarly we define $g(x)$ and denote by 
\begin{equation} \label{defofcn}
c(x)=f^{-1}(g(x)) \, .
\end{equation}

\begin{Proposition}\label{Dj}
Let $\big ( d_{i,n} )_{ 1 \leq i \leq n}$ be a triangular  array of real-valued square-integrable martingale differences, adapted to a triangular array of filtrations $({\mathcal F}_{i,n})_{ 0 \leq i \leq n}$. 
Let $(b_n)_{n\geq 0}$ be a sequence of positive numbers  satisfying \eqref{condseq}, and let
$$
\bar Z_{n}=
\left \{ \frac{d_{1,n}+ \cdots + d_{[nt],n}}{b_n}, t \in [0,1] \right \} \, .
$$
Assume that  the  three following conditions holds
\begin{enumerate}
\item There exists a positive number $V$ such that, for any $\delta>0$ and any $t \in [0,1]$, 
\begin{equation} \label{condvarMDPFlo}
\limsup_{n \rightarrow \infty} \frac{n}{b_n^2} \log {\mathbb P}  \left (  \left| \left ( \frac 1 n  \sum_{i=1}^{[nt]} {\mathbb E}(d_{i,n}^2|{\mathcal F}_{i-1,n}) \right ) - Vt \right |> \delta \right ) = -\infty \, .
\end{equation}
\item  For any $\varepsilon>0$ and $\delta >0$, 
\begin{equation} \label{lindMDPFlo}
\limsup_{n \rightarrow \infty} \frac{n}{b_n^2} \log {\mathbb P} \left ( \frac 1 n  \sum_{i=1}^{n} {\mathbb E}\left (d_{i,n}^2 {\bf 1}_{|d_{i,n}| > \varepsilon  n b_n^{-1}}\Big | {\mathcal F}_{i-1,n} \right) > \delta \right ) = -\infty \, .
\end{equation}
\item 
\begin{equation} \label{onecondm}
\frac{n}{b_n^2} \log \left  (\sup_{n\leq m\leq c(n+1)}\sup_{1\leq k\leq m } n  \big \Vert  {\mathbb{P}}%
\big (|d_{k,m}|>b_n  | {\mathcal F}_{k-1,m} \big )  \big   \Vert_{\infty} \right )  \rightarrow -\infty \text{ as }n\rightarrow
\infty \,,  
\end{equation}
where $c(n)$ is defined in \eqref{defofcn}.
\end{enumerate}
Then, for any 
Borel set $\Gamma \subset D([0,1])$,
\begin{equation}\label{fmdpbis}
- \inf_{ \varphi  \in \Gamma^\circ}  I_V(\varphi) \leq \liminf_{n \rightarrow \infty} 
\frac{n}{b_n^2} \log 
{\mathbb P} \left (
 \bar Z_{n} \in \Gamma 
\right ) 
\leq 
\limsup_{n \rightarrow \infty} 
\frac{n}{b_n^2}  \log 
{\mathbb P} \left ( \bar Z_{n}  \in \Gamma 
\right )  \leq - \inf_{ \varphi  \in \bar \Gamma}  I_V(\varphi) \, ,
\end{equation}
where $I_V$ is defined as in Theorem \ref{MDPth}.
\end{Proposition}

Let us conclude the proof of Theorem \ref{MDPth}. Let $(x_n)_{n \geq 1}$ be any sequence of points in $S^{d-1}$. We apply Proposition \ref{Dj} to the   
martingale differences 
$d_{i,n}=D_{i, x_n}$ (recall that $D_{i,x}$ is the martingale difference of the decomposition  \eqref{MartCob1}).  Condition \eqref{lindMDPFlo} is clearly satisfied thanks
to \eqref{MartCob1} and the fact that 
\begin{multline*}
\left \| {\mathbb E} \left ( X_{k,x_n}^2 {\bf 1}_{|X_{k,x_n}| > \varepsilon  n b_n^{-1}} \Big | 
\mathcal{F}_{k-1} \right ) \right \|_\infty 
=
\left \| \int  (\sigma(g, A_{k-1}x_n))^2  {\bf 1}_{|\sigma(g, A_{k-1}x_n)| > \varepsilon  n b_n^{-1}} \mu (dg) 
\right \|_\infty \\
\leq 
\int_G  (\log N(g))^2  {\bf 1}_{\log N(g)> \varepsilon  n b_n^{-1}}\mu(dg)  \, .
\end{multline*}

To check Condition \eqref{condvarMDPFlo}, we apply Proposition 3.1 in \cite{BQ}, which implies that, for any $\delta>0$ and any $t \in [0,1]$, there exist $A>0$ and  $\alpha>0$ such that, 
for the variance $V$ defined in Theorem  \ref{MDPth},
$$
{\mathbb P}  \left (  \left| \left ( \frac 1 n  \sum_{i=1}^{[nt]} {\mathbb E}(D_{i,x_n}^2|{\mathcal F}_{i-1}) \right ) - Vt \right |> \delta \right ) \leq A e^{-\alpha n} \, .
$$
Condition \eqref{condvarMDPFlo} follows then easily, since $n^2 b_n^{-2} \rightarrow \infty$ as $n \rightarrow \infty$. 

It remains to check Condition \eqref{onecondm}. By  \eqref{MartCob1} again, it is equivalent to prove the condition for $X_{k,x_m}$ instead of $D_{k,x_m}$. 
Now 
\begin{equation*}\label{boundX0}
\left \| {\mathbb E} \left ( {\bf 1}_{|X_{k,x_m} | \geq b_n} \Big |
\mathcal{F}_{k-1} \right ) \right \|_\infty 
=
\left \|  \int {\bf 1}_{\sigma(g, A_{k-1}x_m) \geq b_n} \, \mu (dg) 
\right \|_\infty 
\leq 
\mu \left \{ \log N \geq  b_n \right \}  \, ,
\end{equation*}
and the result follows by \eqref{ARC}.  $\Diamond$  

\section{The case of weak moment of order $p>1$}  
In this section,  we study the asymptotic behavior of  \eqref{modev} 
when $\log N(Y_1)$ has only  a weak moment of order $p>1$. 
\label{weakMom}
\setcounter{equation}{0}
\begin{Theorem}\label{devweak}
Let $p > 1$ and and assume that
\begin{equation}\label{weakM}
 \sup_{t>0}  t^p \mu \left \{ \log N >t \right \} <\infty \, .
\end{equation}
Let $\alpha \in (1/2,1]$ and $\alpha\geq 1/p$. 
Then there exists a positive constant $C$ such that, 
for  any $y>0$,
\begin{equation}\label{moderate}
\limsup_{n \rightarrow  \infty} 
n^{\alpha p-1} 
\sup_{\|x\|=1}
{\mathbb P} \
\left ( \max_{1 \leq k \leq n} | \log \| A_k x \|-k\lambda_\mu |>
n^\alpha y \right ) \leq \frac{C}{y^p} \, .
\end{equation}
\end{Theorem}
\noindent{\bf Proof of Theorem \ref{devweak}.} 

{\it  The case $p>2$.} In that case 
the decomposition \eqref{MartCob2} holds, and 
it is equivalent to 
prove \eqref{moderate}  for $M_{k,x}$ instead of 
$(\log \| A_k x \|-k \lambda_\mu)$. To do this, we shall 
apply the following inequality due to Haeusler \cite{H}:
for all $\gamma ,u,v>0$,
\begin{multline}\label{Fuk}
{\mathbb P}\left ({\max_{1 \leq k \leq n} \left | M_{k,x}
\right | \geq \gamma } \right )\leq \sum_{i=1}^n 
{\mathbb P} \left (| D_{i,x}| \geq u \right ) 
+ 2 {\mathbb P} \left (\sum_{i=1}^n {\mathbb E} (D_{i,x}^2| 
{\mathcal F}_{i-1} ) \geq v \right ) \\ + 
2 \exp \left(\gamma u^{-1}
\left  (1 - \log \left ( \gamma u v^{-1} \right ) \right )\right ) \, .
\end{multline}
Note that if \eqref{weakM} holds for $p>2$, then 
\begin{equation}\label{boundX}
\left \| {\mathbb E} \left ( X_{k,x}^2 | 
\mathcal{F}_{k-1} \right ) \right \|_\infty 
=
\left \| \int  (\sigma(g, A_{k-1}x))^2 \mu (dg) 
\right \|_\infty 
\leq 
\int_G  (\log N(g))^2 \mu(dg) < \infty \, ,
\end{equation}
and there exists a positive constant $C$ such that
\begin{equation}\label{boundX2}
\left \| {\mathbb E} \left ( {\bf 1}_{|X_{k,x} | \geq u} \Big |
\mathcal{F}_{k-1} \right ) \right \|_\infty 
=
\left \|  \int {\bf 1}_{\sigma(g, A_{k-1}x) \geq u} \, \mu (dg) 
\right \|_\infty 
\leq 
\mu \left \{ \log N \geq  u \right \} \leq \frac{C}{u^p} \, ,
\end{equation}
for any $u>0$ and any positive integer $k$.
Using again that $\psi$ is bounded  we infer from \eqref{MartCob1} and \eqref{boundX}  that there
 exist two positive  constants $c_1, c_2$ such that
\begin{align}\label{boundM2}
\sup_{\|x\|=1} 
\left \| {\mathbb E} \left (  D_{k,x}^2  | 
\mathcal{F}_{k-1} \right ) \right \|_\infty &\leq c_1 \, ,
\\
\sup_{u>0} u^p 
\sup_{\|x\|=1}  {\mathbb P} \left (| D_{k,x}| \geq u \right ) 
&\leq 
\sup_{u>0} u^p 
\sup_{\|x\|=1} 
\left \| {\mathbb E} \left ( {\bf 1}_{|D_{k,x}| \geq u}| 
\mathcal{F}_{k-1}\right ) \right \|_\infty
 \leq c_2 \, , \label{boundM3}
\end{align}
for any positive integer $k$. 
Taking 
$\gamma=n^\alpha y 
$,
  $u=  n^\alpha y/r$ with $r \in (0, \infty)$, and $
v= 2n c_1$
in \eqref{Fuk},
we get
\begin{equation}\label{weakLM}
{\mathbb P}\left ({\max_{1 \leq k \leq n} \left | M_{k,x}
\right | \geq  n^\alpha y } \right )\leq c_3 \left ( 
\frac{1}{y^p n^{\alpha p
-1}} + \frac{1}{y^{2r} n^{(2\alpha -1) r}}\right) \, ,
\end{equation}
for some positive constant $c_3$. 
Selecting $r > (\alpha p -1)/(2 \alpha -1)$, the upper bound \eqref{moderate} follows directly from
\eqref{weakLM}.

\medskip

{\it The case $p \in (1,2)$.} Let $x\in S^{d-1}$. We have 
$$
\log \|A_n x\| -n\lambda_\mu=\sum_{k=1}^n (X_{k, x}-\lambda_\mu)
=\sum_{k=1}^n (D_{k, x}+R_{k, x})
:=M_{n,x}+U_{n,x}\, ,$$
where
$$
D_{n,x}= \sigma(Y_n, A_{n-1}x)-\int_G\sigma(g, A_{n-1}x)\mu(dg), \ \text{and} \
R_{n,x}=\int_G\sigma(g, A_{n-1}x)\mu(dg)-\lambda_\mu \, .
$$
Notice  that ${\mathbb E}(D_{k,x}|{\mathcal F}_{k-1})=0$. We use now the basic inequality
\begin{equation}\label{trivial}
\P\left (\max_{1\le k\le n}|\log \|A_k x\| -k\lambda_\mu|\ge n^\alpha y \right) \le 
\P \left (\max_{1\le k\le n}|M_{k,x}|\ge n^\alpha y/2\right )
+\P\left (\max_{1\le k\le n}|U_{k,x}|\ge n^\alpha y/2 \right )\, . 
\end{equation}
We first deal with the second term on the right-hand side of \eqref{trivial}. 
We shall need the following extension of  Theorem 3 in 
\cite{WZ}.  
The proof is given in  Appendix. 
\begin{Theorem}
\label{LpmainNSpmoinsque2} Let $p \in ]1, 2[ $ and  $(X_k)_{k \in {\mathbb Z}}$ be a sequence of real-valued random variables in ${\mathbb L}^p $ and adapted to a non-decreasing filtration $(\F_k)_{k \in {\mathbb Z}}$. Let 
$S_i = X_1 + \cdots + X_i$ and 
$S_n^* = \max_{1 \leq i \leq n} |S_i|$.  Then, for any $n \geq 1$, 
\begin{equation} \label{eq:bs0314NSpmoinsque2}
\Vert S_n^* \Vert_p \leq  ( 2 c_p  + 1) \left (  \sum_{j=1}^{n} \Vert X_j \Vert_p^p \right )^{1/p} + 2^{(p-1)/p} ( 2 c_p  + 1)  \sum_{j=0}^{r-1}  
\left ( \sum_{k=1}^{2^{r-j}} \Vert {\mathbb E} ( S_{k 2^j } -S_{(k-1) 2^j }  | {\mathcal F}_{(k-1) 2^j} ) \Vert_p^p \right )^{1/p} \, ,
\end{equation}
where $c_p = 2^{1/p} \frac{p}{p-1}$ and $r$ is the unique positive integer such that $2^{r-1} \leq n < 2^r$.  
\end{Theorem}

For $k\le 0$, set $R_{k,x}=R_{0,x}$ and $\F_k=\F_0$. Observe that  $|R_{k,x}|\le \int_G\log N(g)  \mu(dg) <\infty$ for every $k\ge 0$. Hence we may apply Theorem \ref{LpmainNSpmoinsque2}  with $X_k:= R_{k,x}$. With that choice, we have 
$$
S_{k 2^j } -S_{(k-1) 2^j } =\sum_{\ell =1}^{2^j} R_{(k-1) 2^j+\ell,x}\, ,
$$
and, using independence (twice), 
$$
\left|\E( R_{(k-1) 2^j+\ell,x}|{\mathcal F}_{(k-1) 2^j})\right|\le 
\sup_{y\in S^{d-1}}\left| \int_G \left(\E\left(\sigma(g, A_{\ell-1} y)\right)
-\lambda_\mu\right) \mu (dg) \right|= 
\, \sup_{y\in S^{d-1}} \left|\E(X_{\ell,y})-\lambda_\mu\right|.
$$
Let $n\ge 1$ and $r\ge 1$ be such that $2^{r-1}\le n< 2^r$. We infer that there exists $C_p>0$, such that 
\begin{align}
\nonumber \left \Vert\max_{1\le k\le n}|U_{k,x}|\, 
\right \Vert_p  &  \le C_p n^{1/p}+ C_p\sum_{j=0}^{r-1} 
2^{(r-j)/p} \sum_{\ell =1}^{2^j}\sup_{\|y\|=1}\left|\E(X_{\ell,y})-\lambda_\mu\right|\\
 \label{majH} &  \le C_{p}n^{1/p} +\frac{C_p2^{1/p}}{2^{1/p}-1} n^{1/p}\sum_{\ell \ge 1} 
  \frac{\sup_{\|y\|=1}\left|\E(X_{\ell,y})-\lambda_\mu\right|}{\ell^{1/p}}
  \, .
\end{align}
Recall that  \eqref{weakM} holds for $p \in (1,2)$. Hence, for any $r<p$, by (6) of \cite{CDJ}, 
\begin{equation}\label{mainbound}
\sum_{n\ge 1} n^{r-2} \sup_{ \|y\|=1}\left |\E (X_{n,y})-\lambda_\mu \right |<\infty\, .
\end{equation}
Since $p+1/p>2$ one can choose $r$ close enough to $p$ in such a way that
$-1/p<r-2$. In particular, it follows that
\begin{equation}\label{mainbound2}
\sum_{\ell\ge 1} \ell^{-1/p} \sup_{\|y\|=1}\left |\E (X_{\ell,y})- \lambda_\mu \right |<\infty\, .
\end{equation}
Hence, using \eqref{majH},
\begin{align*}
\left \|\max_{1\le i\le n} |U_{i,x}|\right \|_p 
  \le \tilde C_p n^{1/p}\, ,
\end{align*}
and 
\begin{equation}\label{secondterm}
\sup_{\|x\|=1}\P\left (\max_{1\le k\le n}|U_{k,x}|\ge y n^{\alpha}/2 \right )
\le \frac{(2\tilde C_p)^p n}{y^pn^{p\alpha}}\, ,
\end{equation}
which ends the control of the second term on the right-hand side of  \eqref{trivial}. 

We now deal with the  first term on the right-hand side of \eqref{trivial}, that is the martingale term. We shall 
need the following result (to be proved in Appendix). It
is a maximal-version of Theorem 2.5 in \cite{GM} (a von Bahr-Esseen inequality
for martingales having  weak moments of order $p \in (1,2)$). 
For a real-valued random variable $X$, let 
$\|X\|_{p,\infty}=\sup_{t>0} t ({\mathbb P}(|X|>t))^{1/p}$. 

\begin{Proposition}\label{weakvbe}
Let $(D_n)_{n\in {\mathbb N}}$ be a sequence of  $({\mathcal F}_n)_{n\in 
{\mathbb N} }$-martingale differences in weak-$L^p$, $p \in (1,2)$. Then 
$$
\P\left (\max_{1\le k\le n} \left |\sum_{j=1}^k  D_j
\right |\ge y \right ) 
\le \frac{K }{y^p}\sum_{k=1}^n \|D_k\|_{p,\infty}^p\, ,
$$
where $K=4p/(p-1) + 8/(2-p)$. 
\end{Proposition}

Now, since \eqref{weakM} holds, then so does \eqref{boundM3}. It follows 
from Proposition \ref{weakvbe} that 
\begin{equation}\label{firstterm}
\sup_{\|x\|=1}\P\left (\max_{1\le k\le 2^n}|M_{k,x}|\ge n^{\alpha}y /2 \right )
\le \frac{C}{y^p n^{\alpha p -1}}\, ,
\end{equation}
for some positive constant $C$. The upper bound \eqref{moderate} follows from \eqref{trivial}, \eqref{secondterm}
and \eqref{firstterm}.

\medskip

{\it The case $p=2$.} We start  from \eqref{trivial}.
Note first that the upper bound \eqref{secondterm}
 still holds for $p=2$, with the same proof. We now deal with the 
 first term on the right-hand side of \eqref{trivial}. Instead 
of Proposition \eqref{weakvbe}, we shall use the following result
of Hao and Liu \cite{HL} (see also Theorem 14 in \cite{CDJ}):
if ${\mathbb P}(|D_{k,x}|>y) \leq {\mathbb P}(X>y)$ for any $y>0$ and some 
positive random variable $X$, then, 
for every
$q>1$, every $\gamma\in (1,2]$ and every $L\in {\mathbb N}$, there exists  $C>0$,  such that for every $n\ge 1$ and every $\lambda >0$, 
\begin{multline}
 \label{complete-conv}\P\left (\max_{1\le k\le n} |M_{k,x}|\ge \lambda \right) 
 \le n \P\left(X>   \frac{\lambda}{4(L+1)}\right)
 \\  +\frac{C}{(\lambda)^{q\gamma(L+1)/(q+L)}} 
\left  \| \E(|D_{1,x}|^\gamma|\F_0)+\cdots + \E(|D_{n,x}|^\gamma|\F_{n-1}) \right \|_q ^{q(L+1)/(q+L)}\, .
\end{multline}
We  apply \eqref{complete-conv} with
$X=\log N(Y_1)+{\mathbb E}(\log N(Y_1))$. 
Since \eqref{weakM} holds with $p=2$, then $X$ has a weak   moment of order $2$, and,  for every $ \gamma \in (1,2)$, there exists $C_{\gamma}>0$ such that for every $n\ge 1$, 
 $$
 \|\E(|D_{n,x}|^\gamma|\F_{n-1})\|_\infty \le C_{\gamma}\, .
 $$
Hence, for every integer $L$ and every $q>1$, there exists 
 $C>0$ such that 
 $$
\P\left (\max_{1\le k\le n}|M_{k,x}|\ge  n^\alpha y /2 \right ) 
\le n \P\left(X\ge \frac{n^\alpha y}{8(L+1)}\right) 
+\frac{C n^{q(L+1)/(q+L)}}{(n^\alpha y)^{ q\gamma (L+1)/(q+L)}}\, .
 $$
 Since $\alpha>1/2$ one may find $\gamma \in (1,2)$, such that 
 $\gamma\alpha >1$. For such a choice, taking $q=L$ large enough, 
 we obtain the desired result.  $\Diamond$

\section{The case of strong moments of order $p\geq 1$}
In this section, we prove some results in the spirit of Baum and Katz \cite{BK} for  the quantity \eqref{modev}. 
\label{strongMom}
\setcounter{equation}{0}
\begin{Theorem}\label{devstrong}
Let $p \geq 1$ and  assume that
\begin{equation}\label{strongM}
 \int (\log N(g))^p \mu (dg) <\infty \, .
\end{equation}
Let $\alpha \in (1/2,1]$ and $\alpha\geq 1/p$. 
Then 
for  any $y>0$
\begin{equation}\label{moderate2}
\sum_{n \geq 1} 
n^{\alpha p-2} 
\sup_{\|x\|=1}
{\mathbb P} \
\left ( \max_{1 \leq k \leq n} | \log \| A_k x \|-k\lambda_\mu |>
n^\alpha y \right ) < \infty \, .
\end{equation}
\end{Theorem}
\begin{Remark}
Theorem \ref{devstrong} is due to Benoist-Quint \cite{BQ} in the case where 
$\alpha=1$ and $p>1$.
\end{Remark}
\begin{Remark}
Let us recall a well known consequence of \eqref{moderate2}, when  $p\in [1,2)$ and $\alpha=1/p$. 
The sequence $\max_{1 \leq k \leq n} | \log \| A_k x \|-k\lambda_\mu |$ being non-decreasing, Inequality 
\eqref{moderate2} with $\alpha=1/p$ is equivalent to 
\begin{equation}\label{moderate2bis}
\sum_{N \geq 1} 
\sup_{\|x\|=1}
{\mathbb P} \
\left ( \max_{1 \leq k \leq 2^N} 
| \log \| A_k x \|-k\lambda_\mu |>
2^{N/p} y \right ) < \infty \, .
\end{equation}
This implies  that, for any $x \in S^{d-1}$, the sequence$ (2^{-N/p}\max_{1 \leq k \leq 2^N} | \log \| A_k x \|-k\lambda_\mu |)_{N \geq 1}$
converges completely. It follows easily that, for any $x \in S^{d-1}$, $n^{-1/p} ( \log \| A_k x \|-k\lambda_\mu)$ converges to $0$ almost surely as $n \rightarrow \infty$. 
Hence  \eqref{moderate2} is a more precise statement than Theorem 1(i) of \cite{CDJ}.
\end{Remark}

Of course, \eqref{moderate2} does not hold for $p=2$ and $\alpha=1/2$. Instead, we have the following result, which implies a bounded law of the iterated logarithm. 

\begin{Theorem}\label{devLIL}
Assume that
$
 \int (\log N(g))^2 \mu (dg) <\infty 
$,
and let $V$ be defined as in Theorem \ref{MDPth}. 
Then 
for  any $y> \sqrt{V} \ge 0$, we have
\begin{equation}\label{BLIL}
\sum_{n \geq 1} 
\frac 1 n 
\sup_{\|x\|=1}
{\mathbb P} \
\left ( \max_{1 \leq k \leq n} | \log \| A_k x \|-k\lambda_\mu |>
 y\sqrt {2n\log \log n} \right ) < \infty \, .
\end{equation}
\end{Theorem}

\begin{Remark}
From \eqref{BLIL} one can easily infer that, for any $x \in S^{d-1}$, 
$$
  \limsup_{n \rightarrow \infty} \frac{|\log \|A_nx\|-n\lambda_\mu|} {\sqrt {2n\log \log n}} \leq  \sqrt{V}, \quad  \text{almost surely.}
$$
Of course, this is a less precise result than the compact law of the iterated logarithm, which also holds 
provided $
 \int (\log N(g))^2 \mu (dg) <\infty 
$ (for instance, this 
is a consequence of Theorem 1(iii), case $p=2$, of \cite{CDJ}). Note however that \eqref{BLIL} and the compact law of the iterated logarithm are 
two different results, which cannot be deduced  from one another. 
\end{Remark}

\noindent{\bf Proof of Theorem \ref{devstrong}.} The proof follows the line of that of Theorem  \ref{devweak}. 

{\it The case $p\geq 2$.} In that case 
the decomposition \eqref{MartCob2} holds, and 
it is equivalent to 
prove \eqref{moderate2}  for $M_{k,x}$ instead of 
$(\log \| A_k x \|-k \lambda_\mu)$. 

Starting from \eqref{MartCob1} and \eqref{boundX2}, we see that 
\begin{equation}\label{stocbound}
\sup_{\|x\|=1, k\geq 1}  {\mathbb P} \left (| D_{k,x}| \geq 
n^\alpha y/r \right )
 \leq \mu \left \{ \log N \geq  n^\alpha y/2r \right \} +
 {\bf 1}_{n^\alpha y \leq 2r (2 \|\psi\|_\infty + |\lambda_\mu |)} \, .
\end{equation}
Taking 
$\gamma=n^\alpha y 
$,
  $u=  n^\alpha y/r$ with $r>0$, and $
v= 2 n c_1$ (cf. \eqref{boundM2} for the definition of  $c_1$)
in \eqref{Fuk},
we get
\begin{equation}\label{strongLM}
{\mathbb P}\left ({\max_{1 \leq k \leq n} \left | M_{k,x}
\right | \geq  n^\alpha y } \right )\leq 
n \mu \left \{ \log N \geq  n^\alpha y/2r \right \} +
n  {\bf 1}_{n^\alpha y \leq  2r (2 \|\psi\|_\infty + |\lambda_\mu |)}
+ \frac{\kappa_1}{y^{2r} n^{(2\alpha -1) r}} \, ,
\end{equation}
for some $\kappa_1>0$. Interverting the sum and the integral, 
we see that 
\begin{equation}\label{momentN}
\sum_{n>0} n^{\alpha p -1} 
\mu \left \{ \log N \geq  n^\alpha y/2r \right \}  \leq 
\frac{\kappa_2}{y^p} \int (\log N(g))^p \mu (dg) \, ,
\end{equation}
for some positive constant $\kappa_2$ depending 
only on $r$. Taking $r > (\alpha p -1)/(2 \alpha -1)$ in 
\eqref{strongLM} and using the upper bound 
\eqref{momentN}, the proof of \eqref{moderate2}
is complete for $p\geq 2$.

\medskip

{\it The case $p \in (1,2)$.} We start again from \eqref{trivial}.
If \eqref{strongM} holds for $p \in (1,2)$ then, by (6) of \cite{CDJ}, \eqref{mainbound}
holds with $r=p$. Since $p+1/p>2$ and $p<2$, there exists $q$ such that $p<q<2$ 
and $ p+1/q>2$. 
Hence, we infer that
\begin{equation}\label{mainbound2bis}
\sum_{n\ge 1} n^{-1/q} \sup_{\|x\|=1}\left |\E (X_{n})- \lambda_\mu \right |<\infty\, , 
\end{equation}
and, using \eqref{majH} in ${\mathbb L}^q$ rather than in ${\mathbb L}^p$, we infer that 
for every $n\ge 1$
\begin{align*}
\left \|\max_{1\le i\le n} |U_{i,x}|\right \|_q 
  \le C_q n^{1/q}\, ,
\end{align*}
for some $C_q>0$.
Hence, 
$$
\sup_{\|x\|=1}\P\left (\max_{1\le k\le n}
|U_{k,x}| \geq y n^{\alpha }/2
\right )\le \frac{(2C_q)^q n}{y^qn^{q\alpha}}\, ,
$$
and, since $q>p$, 
\begin{equation}
\sum_{n \geq 1} 
n^{\alpha p-2} 
\sup_{\|x\|=1}
{\mathbb P} \
\left ( \max_{1 \leq k \leq n} | U_{k,x} |>
n^\alpha y/2 \right ) < \infty \, .
\end{equation}

It remains to deal with the first term on the right-hand side of
\eqref{trivial}. Applying  \eqref{MartCob1} and \eqref{boundX2} again, we see that $|D_{k,x}|$ is uniformly (with respect to $k$ and $x$) stochastically bounded by
the random variable $\log N(Y_1) + | \lambda_\mu | +2 \| \psi \|_\infty$.  
Following exactly the proof of Theorem 2 of \cite{DM}, it follows that 
\begin{equation}
\sum_{n \geq 1} 
n^{\alpha p-2} 
\sup_{\|x\|=1}
{\mathbb P} \
\left ( \max_{1 \leq k \leq n} |  M_{k,x} |>
n^\alpha y/2 \right ) < \infty \, ,
\end{equation}
and the proof is complete for $p \in (1,2)$.

\medskip

{\it The case $p=1$.} In that case $\alpha=1$ also. So, let us start from the inequality \eqref{trivial} with $\alpha=1$.

The second term on the right-hand side of \eqref{trivial} may be handled thanks to a ``maximal version" of Proposition 3.1 of \cite{BQ}, which implies that: for any $y>0$, there exists
constants $A>0$, $\alpha>0$ such that 
\begin{equation}\label{maxBQ}
\sup_{\|x\|=1} \P\left (\max_{1\le k\le n}|U_{k,x}|\ge n y/2 \right ) \leq A e^{-\alpha n} \, .
\end{equation}
Note that a direct application of of Proposition 3.1 of \cite{BQ}
gives us \eqref{maxBQ} without the maximum over $k$. 
To prove  \eqref{maxBQ},
it is convenient  to work on the projective space $X:=\P_{d-1}(\R)$ rather than 
on $S^{d-1}$. Denote by $\bar x$ the class (in $X$) of $x\in {\mathbb R}^d- \{0\}$. 
Then, we define $U_{n,\bar x}:=U_{n,x}$. Recall that, by our assumptions, 
the Markov chain $(\overline{A_{n-1}\cdot  x})_{n\ge 1}$ with (compact) state  space 
$X$ and transition probability given by $Pf(\bar x):= 
\int_Gf(\overline{g\cdot  x})\mu(dg)$ has a unique invariant probability $\nu$. 
In particular, for every continuous function  $f$ on $X$,  the sequence $$\left (\frac 1 n \E \left (\sum_{k=1}^nf \left (\overline{A_{k-1}\cdot 
 x} \right ) \right ) \right )_{n\ge 1}$$ converges uniformly (with respect to $\bar x$) 
to $\nu(f)$. 
Hence, there exists an integer $m\ge 1$ such that 
\begin{equation}\label{unif}
\sup_{\bar x\in X} |\E(U_{m,\bar x})|< my/4\, .
\end{equation}
The maximal inequality \eqref{maxBQ} follows then by applying Lemma 23  in \cite{DMPU}.

 Let us deal with the first term on the right-hand side of \eqref{trivial}. 
Let $\Gamma_n:= \cup_{k=1}^n \{ \log(N(Y_k))\ge yn\}$ and notice that 
on $\Gamma_n^c$, $$\sigma(Y_k, A_{k-1}x)=\sigma(Y_k, A_{k-1}x){\bf 1}_{\{
\log(N(Y_k))< yn\}}\, .$$ 
Define 
$$
\tilde M_{k,x}:= \sum_{j=1}^k \Big(\sigma(Y_k, A_{k-1}x){\bf 1}_{\{
\log(N(Y_k))< yn\}}- \int_G\sigma(g, A_{k-1}x){\bf 1}_{\{
\log(N(g))< yn\}}\mu(dg) \Big)\, ,
$$
and note that 
$(\tilde M_{k,x})_{1\le k\le n}$ is a martingale.
Let 
$$
I(n)=\int_G \log (N(g)){\bf 1}_{\{N(g)\ge yn\}} \mu(dg)\, ,
$$
and note that
$$
\P \left (\max_{1\le k\le n}|M_{k,x}|\ge n y/2\right )
\le \P(\Gamma_n)+ \P \left (\max_{1\le k\le n}|\tilde M_{k,x}|\ge n y/4\right ) 
+ {\bf 1}_{\{I(n) \geq y/4\}} \, .
$$
Using Doob's maximal inequality, we infer that 
$$
\P \left (\max_{1\le k\le n}|M_{k,x}|\ge n y/2\right )
\le  \P(\Gamma_n)+ \frac{64}{n^2 y^2} \E\left (\tilde M_{n,x}^2 \right )
+   {\bf 1}_{\{I(n) \geq y/4\}} \, .
$$
The last term on the right-hand side  is equal to 0 for $n$ large enough since \eqref{strongM} holds 
with $p=1$. 

Now, $\P(\Gamma_n)\le n\P(\log N(Y_1)\ge yn)$. Hence it is standard that 
$\sum_{n\ge 1}n^{-1} \P(\Gamma_n)<\infty$, since \eqref{strongM} holds 
with $p=1$.

On the other hand, 
$$
\sup_{\|x\|=1}  \E\left (\tilde M_{n,x}^2 \right ) \le n\int_G (\log N(g))^2{\bf 1}_{\{\log(N(g))< yn\}} 
\mu(dg)\, .
$$
Then, it is also standard that $\sum_{n\ge 1} 
n^{-3}\sup_{\|x\|=1}  \E(\tilde M_{n,x}^2)<\infty$, since \eqref{strongM} holds 
with $p=1$. $\Diamond$

\medskip

\noindent{\bf Proof of Theorem \ref{devLIL}.}  In that case 
the decomposition \eqref{MartCob2} holds, and 
it is equivalent to 
prove \eqref{moderate2}  for $M_{k,x}$ instead of 
$(\log \| A_k x \|-k \lambda_\mu)$. Moreover, we have 
$V=\E(D_{1,x}^2)$.

\smallskip
We shall proceed by truncation. Let $y> \sqrt{V}$ and set $\varepsilon 
:=y-\sqrt{V}$.
Let $n\ge 1$. Let $\alpha>0$ be fixed for the moment. For every $
1\le k \le n$, define 
$$
\tilde D_{k,n, x}:=D_{k,x}{\bf 1}_{\{|D_{k,x}|\le \alpha \sqrt n/\sqrt {\log\log n}\}}
-\E\left (D_{k,x}{\bf 1}_{\{|D_{k,x}|\le \alpha \sqrt n/\sqrt {\log\log n}\}}
|\F_{k-1} \right )
$$ 
and 
$$
\tilde M_{k,n,x}=\sum_{j=1}^k \tilde D_{j,n,x}\, .
$$
Then, using Markov's  inequality and stationarity
\begin{align*}
\P\left(\max_{1\le k \le n}|M_{k,x}|>y\sqrt {2n\log \log n} \right)
\le & \ \P\left(\max_{1\le k \le n}|\tilde M_{k,n,x}|>(y-\varepsilon/2)
\sqrt {2n\log \log n} \right)\\ & +\P\left(\max_{1\le k \le n}|M_{k,x}-\tilde M_{k,n,x}|>y\sqrt {2n\log \log n} \right)\\ 
\le & \ \P\left(\max_{1\le k \le n}|\tilde M_{k,n, x}|>(\varepsilon/2)
\sqrt {2n\log \log n} \right) \\ & +\frac{2n}{\varepsilon
\sqrt {2n\log \log n}}\E \left (2|D_{1,x}|{\bf 1}_{\{|D_{1,x}| > \alpha \sqrt n/\sqrt {\log\log n}\}} \right )\, .
\end{align*}
Now, starting from  \eqref{MartCob1}  and arguing as in   \eqref{boundX2},  
\begin{multline*}
 \sup_{\|x\|=1, k \geq 1}\E\left (|D_{1,x}|{\bf 1}_{\{|D_{1,x}| > \alpha \sqrt n/\sqrt {\log\log n}\}} \right ) \\
\leq \E\left ((N(Y_1) + |\lambda_\mu| + 2 \|\psi\|_\infty) {\bf 1}_{\{N(Y_1) + |\lambda_\mu| + 2 \|\psi\|_\infty > \alpha \sqrt n/\sqrt {\log\log n}\}} \right ) \, .
\end{multline*}
Since $
 \int (\log N(g))^2 \mu (dg) <\infty 
$, it is now standard that 
$$
\sum_{n \geq 1} \frac{4}{\varepsilon
\sqrt {2n\log \log n}} \sup_{\|x\|=1} \E\left (|D_{1,x}|{\bf 1}_{\{|D_{1,x}| > \alpha \sqrt n/\sqrt {\log\log n}\}} \right ) <\infty \, .
$$

Hence, it remains to prove that 
\begin{equation}\label{convseries}
\sum_{n\ge 1} \frac 1 n \sup_{\|x\|=1} \P\left(\max_{1\le k \le n}|\tilde M_{k,n,x}|>
(y-\varepsilon/2)
\sqrt {2n\log \log n} \right) <\infty\, .
\end{equation}

We shall use the following sharper version of Haeusler's bound 
\eqref{Fuk} (see the end of the Proof of Lemma 1 in \cite{H}).
\begin{multline}\label{Fuk2}
{\mathbb P}\left ({\max_{1 \leq k \leq n} \left | \tilde M_{k,n,x}
\right | \geq \gamma } \right )\leq \sum_{i=1}^n 
{\mathbb P} \left (| \tilde D_{i,n,x}| \geq u \right ) 
+ 2 {\mathbb P} \left (\sum_{i=1}^n {\mathbb E} (\tilde D_{i,n,x}^2| 
{\mathcal F}_{i-1} ) \geq v \right ) \\ + 
2 \exp \left(\gamma u^{-1}-(\gamma u^{-1}+vu^{-2})\log(\gamma uv^{-1}+1)
\right ) \, .
\end{multline}
Notice that for every $1\le k\le n$, $|\tilde D_{k,n,x}|\le2 \alpha \sqrt n/\sqrt {\log\log n}$ and that 
\begin{equation}\label{triv-est}
\sum_{i=1}^n {\mathbb E} (\tilde D_{i,n,x}^2| 
{\mathcal F}_{i-1} )\le \sum_{i=1}^n {\mathbb E} (D_{i,x}^2| 
{\mathcal F}_{i-1} )\, .
\end{equation}
By Proposition 3.1 of  \cite{BQ}, 
\begin{equation}\label{est-BQ}
\sum_{n\ge 1}\sup_{\|x\|=1} {\mathbb P} \left (\sum_{i=1}^n {\mathbb E} (\tilde D_{i,n,x}^2| 
{\mathcal F}_{i-1} ) \geq n(\sqrt V +\varepsilon)^2 \right ) <\infty\, .
\end{equation}

We shall  apply \eqref{Fuk2} with $\gamma:=(y-\varepsilon/2)
\sqrt {2n\log \log n}$, $v:= (\sqrt V+\varepsilon)^2n$ and 
$u:= 4 \alpha \sqrt n/\sqrt {\log\log n}$.
Using that for every $t\ge 0$, $\log(1+t)\ge t-t^2/2$, we infer that
\begin{gather*}
\gamma u^{-1}-(\gamma u^{-1}+vu^{-2})\log(\gamma uv^{-1}+1)
\le -\frac{\gamma^2v^{-1}}2(1-\gamma uv^{-1})\, .
\end{gather*}
Since $$\frac{\gamma^2v^{-1}}{2\log\log n}=\frac{(\sqrt V+\varepsilon/2)^2}{(\sqrt V +
\varepsilon/4)^2}>1 \, , $$ and since  $\gamma uv^{-1} = 4\sqrt 2\alpha (\sqrt V+\varepsilon/2)(\sqrt V
+\varepsilon/4)^2 \rightarrow 0$ as $\alpha \rightarrow 0$, we can choose 
$\alpha $ small enough  in such a way  that there exists $\delta>1$ for which
\begin{equation}\label{cond-bertrand}{\rm exp}\left( \gamma u^{-1}-(\gamma u^{-1}+vu^{-2})\log(\gamma uv^{-1}+1)\right)\le 
(\log n)^{-\delta}\, .
\end{equation}
Combining \eqref{Fuk2}, \eqref{triv-est}, \eqref{est-BQ} and 
\eqref{cond-bertrand} we conclude that \eqref{convseries} holds. 
$\Diamond$

\section{Appendix}

\setcounter{equation}{0}

\subsection{Proof of Proposition \ref{Dj}} As we shall see the proposition is a consequence of the following more general result concerning the functional moderate deviation principle 
of an array of martingale differences.

\begin{Theorem}\label{PWDFMDP}
Let $(d_{i,n})_{1 \leq i \leq n}$ be an array of square-integrable martingale differences, adapted to an array of filtrations $({\mathcal F}_{i,n})_{0\leq i \leq n}$. 
Let $(b_n)_{n\geq 0}$ be a sequence of positive numbers  such that $b_n /{ \sqrt n} \rightarrow \infty$ and $b_n /n \rightarrow 0$ as $n \rightarrow \infty$, and  let
$$
\bar Z_{n}=
\left \{ \frac{d_{1,n}+ \cdots + d_{[nt],n}}{b_n}, t \in [0,1] \right \} \, .
$$
Suppose the  conditions \eqref{condvarMDPFlo} and \eqref{lindMDPFlo} satisfied. In addition, assume that
\begin{equation} \label{negcond1Flo}
\limsup_{n \rightarrow \infty} \frac{n}{b_n^2} \log {\mathbb P} \left ( \max_{ 1 \leq k \leq n }  |d_{k,n}| >  b_n \right ) = - \infty \, .
\end{equation}
and that,  for any  $\lambda>0, \delta>0$, 
\begin{equation} \label{negcond2Flo}
\limsup_{n \rightarrow \infty} \frac{n}{b_n^2} \log {\mathbb P} 
\left ( \frac{n}{b_n^2} \sum_{k=1}^n  {\mathbb E} \left (   e^{ \frac{\lambda b_n | d_{k,n}|}{n } }  {\bf 1}_{n b_n^{-1} < |d_{k,n}| \leq b_n } | {\mathcal F}_{k-1,n} \right ) > \delta \right ) = - \infty \, .
\end{equation}
Then,  the functional moderate deviation principle \eqref{fmdpbis} holds.
\end{Theorem}
\noindent {\bf Proof of Theorem \ref{PWDFMDP}.} The proof will be done with the help of a truncature argument, using Puhalskii's functional moderate deviation principle for the main part and proving that the other parts 
have negligible contributions. 

First, to soothe the notations, we suppress the index $n$ and we denote $d_k = d_{k,n}$ and $\F_k = \F_{k,n}$. We use a tuncation of the variables $d_k$ as follows: for all $1\leq k\leq n$,
let
\begin{equation*}
{\bar d}_{k}:=d_{k} {\bf 1}_{|d_{k}| \leq   n b_n^{-1}}-{\mathbb{E}}%
\left (d_{k} {\bf 1}_{|d_{k}| \leq   n b_n^{-1}} | \F_{k-1} \right ) \,,
\end{equation*}%
\[
d_k':=d_{k} {\bf 1}_{ n b_n^{-1} < |d_{k}| \leq   b_n}-{\mathbb{E}}
\left (d_{k} {\bf 1}_{n b_n^{-1} < |d_{k}| \leq   b_n} | \F_{k-1} \right ) 
\]
and
\begin{equation*}
d_k'':= d_{k} {\bf 1}_{  |d_{k}| >  b_n}-{\mathbb{E}}%
\big (d_{k} {\bf 1}_{ |d_{k}| > b_n} | \F_{k-1} \big )   \,.
\end{equation*}
With these notations, we clearly have that, for any $t \in [0,1]$ 
\[
Z_n (t) = b_n^{-1} \sum_{k=1}^{[nt] } d_k = {\bar Z}_n (t) + Z_n'(t) + Z_n'' (t)   \, , 
\]
with ${\bar Z}_n (t) = b_n^{-1} \sum_{k=1}^{[nt] } {\bar d}_k$, $Z'_n (t) = b_n^{-1} \sum_{k=1}^{[nt] } d_k'$ and  $Z''_n (t) = b_n^{-1} \sum_{k=1}^{[nt] } d_k''$.
Notice first that 
\begin{equation*}
b_n^{-1} \sum_{j=1}^{n }{\mathbb{E}}%
\left (|d_{k} | {\bf 1}_{ |d_{k}| > b_n} | \F_{k-1} \right )  \leq \frac{n}{b_n^{2}} \frac{1}{n} \sum_{j=1}^{n }{\mathbb{E}}%
\left (d^2_{k}  {\bf 1}_{ |d_{k}| > b_n} | \F_{k-1} \right )
\end{equation*}%
and that for any $\delta >0$,
\begin{equation*}
\frac{n}{b_n^2}\log {\mathbb{P}}\left (b_n^{-1} \sum_{k=1}^{n} |d_{k}| {\bf 1}_{  |d_{k}| >  b_n} \geq \delta \right )\leq \frac{n}{b_n^2} \log {\mathbb{P}}\left (\max_{1\leq k\leq
n}| d_{k}|> b_n \right )\,.
\end{equation*}%
Hence, by taking into account conditions \eqref{lindMDPFlo} and \eqref{negcond1Flo}, we can deduce that the process $Z''_n$ has a negligible contribution to the functional moderate deviation principle
(see Theorem 4.2.13 in \cite{DZ}).  

On another hand,  note that 
$({\bar d}_k)_{1\leq k \leq n }$ is a triangular sequence of martingale differences such that $\Vert {\bar d}_k \Vert_{\infty} \leq 2n/b_n$. Using conditions \eqref{condvarMDPFlo} and \eqref{lindMDPFlo}, we can 
apply the functional moderate deviation principle of Puhalskii \cite{Pu} which entails that ${\bar Z}_n$ satisfies \eqref{fmdpbis}.   Therefore to end the proof it remains to show that the process $Z'_n$ has a negligible contribution to the functional 
moderate deviation principle; that is: for any $\delta >0$, 
\begin{equation} \label{negsecondMDP}
\limsup_{n \rightarrow \infty} \frac{n}{b_n^2}  \log 
{\mathbb P} \left ( \sup_{t \in [0,1]} | Z'_{n} (t) | > \delta \right ) = - \infty \, .
\end{equation}
Observe that 
\[
b_n^{-1} \sum_{j=1}^{n }{\mathbb{E}}%
\left (|d_{k} | {\bf 1}_{ n b_n^{-1} < |d_{k}| \leq   b_n} | \F_{k-1} \right )  \leq  \frac{1}{n} \sum_{j=1}^{n }{\mathbb{E}}%
\left (d^2_{k}  {\bf 1}_{ |d_{k}| > nb^{-1}_n} | \F_{k-1} \right ) \, , 
\]
which by using condition \eqref{lindMDPFlo} implies that \eqref{negsecondMDP} will hold if we can prove that, for any $\delta >0$, 
\begin{equation} \label{negsecondMDPp1}
\limsup_{n \rightarrow \infty} \frac{n}{b_n^2}  \log 
{\mathbb P} \left (  b_n^{-1} \sum_{k=1}^n | d_k| {\bf 1}_{ n b_n^{-1} < |d_{k}| \leq   b_n} > \delta \right ) = - \infty \, .
\end{equation}
With this aim, we use the arguments developed  in the proof of Proposition 1 in \cite{Wo}. For the sake of clarity, let us give some details. Take $\lambda $ a positive number and set 
$ Y_{k,\lambda} : =   \frac{2 \lambda b_n}{n}| d_k| {\bf 1}_{ n b_n^{-1} < |d_{k}| \leq   b_n} $.  We have
\begin{multline*}
{\mathbb P} \left (  b_n^{-1} \sum_{k=1}^n | d_k| {\bf 1}_{ n b_n^{-1} < |d_{k}| \leq   b_n} > \delta \right )  = {\mathbb P} \left ( \sum_{k=1}^n  Y_{k,\lambda}   >  \frac{ 2 \lambda \delta b_n^2}{n} \right )  \\
 \leq 
{\mathbb P} \left ( \sum_{k=1}^n  \left \{  Y_{k,\lambda} -
 \log \E \left (   e^{Y_{k,\lambda}} | {\mathcal F}_{k-1} 
 \right )  
\right \}  >  \frac{  \lambda \delta b_n^2}{n} \right )  + 
{\mathbb P} \left ( \sum_{k=1}^n   \log  \E \left (  e^{Y_{k,\lambda}} | {\mathcal F}_{k-1} \right )   >  \frac{  \lambda \delta b_n^2}{n}  \right ) \, .
\end{multline*}
Since the $ Y_{k,\lambda}$ are ${\mathcal F}_k$-measurable, we have \[
\E \left( \frac{\prod_{k=1}^n e^{Y_{k,\lambda}} }{ \prod_{k=1}^n  \E ( e^{Y_{k,\lambda} }| {\cal F}_{k-1}) } \right )= 1 \, .
\]
Hence 
\[
\limsup_{n \rightarrow \infty} \frac{n}{b_n^2}  \log {\mathbb P} \left ( \sum_{k=1}^n  \left \{  Y_{k,\lambda} - \log \E \left (   e^{Y_{k,\lambda}} | {\mathcal F}_{k-1} \right )  \right \}  
>  \frac{  \lambda \delta b_n^2}{n} \right ) \leq - \lambda \delta \, , 
\]
which is going to $-\infty$ by letting $\lambda$ tend to $\infty$. Hence,  to prove \eqref{negsecondMDPp1} (and then \eqref{negsecondMDP}),  it suffices to  show that, for any positive $\lambda$, $\delta$, 
\[
\limsup_{n \rightarrow \infty} \frac{n}{b_n^2}  \log 
{\mathbb P} \left ( \sum_{k=1}^n   \log \E \left (   e^{Y_{k,\lambda}} | {\mathcal F}_{k-1} \right )   >  \frac{  \lambda \delta b_n^2}{n} \right ) = - \infty \, .
\]
This holds under condition \eqref{negcond2Flo} by taking into account that  $e^{x{\bf 1}_A}-1=(e^{x}-1){\bf 1}_A$ and also that $\log (1+u)\leq u$ for any $x>0, u>0$.  The 
proof of Theorem \ref{PWDFMDP} is therefore complete. $\diamond$

\medskip

\noindent {\bf End of the proof of Proposition \ref{Dj}.} We start with some observations. Obviously condition \eqref{negcond2Flo} holds under the stronger one: for any $\lambda >0$,
\[
\limsup_{n \rightarrow \infty} \frac{n}{b_n^2} \sum_{k=1}^n 
 \left  \Vert  {\mathbb E} \left (   e^{ \frac{\lambda b_n | d_{k,n}|}{n } }  {\bf 1}_{n b_n^{-1} < |d_{k,n}| \leq b_n } | {\mathcal F}_{k-1,n}  \right ) \right \Vert_{\infty} =0 \, .
\]
Note now that this condition is equivalent to the following one. There is
a constant $C$ with the following property: for any $\lambda >0$ there exists a positive integer $
N(\lambda )$ such that for $n>N(\lambda )$,
\begin{equation}  \label{negcond2Flostrong}
\frac{n}{b_n^2} \sum_{k=1}^n  \left \Vert {\mathbb{P}} \left (| d_{k,n}| >u n b_n^{-1} | {\mathcal F}_{k-1,n}  \right )  \right \Vert_{\infty} \leq
C \exp (-\lambda u) \quad \text{ for all }1\leq u\leq b_n^2 /n\,.  
\end{equation}
(The proof of this equivalence can be done by following the proof of Comment 6 in \cite{MP}).  To end the proof of the proposition, it remains to show that condition \eqref{onecondm} implies \eqref{negcond2Flostrong} 
(since it obviously implies condition \eqref{negcond1Flo}). Under the regularity conditions \eqref{condseq}, this can be achieved by following the 
lines of the proof of Corollary 7 in \cite{MP} (by taking $s_n= \sqrt n$, $k_n = n$ and $a_n = n/b_n^2$).  $\diamond$

\subsection {Proof of Theorem \ref{LpmainNSpmoinsque2}} 
{\bf Proof.} The proof follows the lines of the proof of Theorem 3 in \cite{WZ}, but in the non-stationary setting, and is then done by induction. For $n=1$, the inequality is clearly true. 
Assume that the inequality holds up to $n-1$ for any sequence $(X_k)_{k \in {\mathbb Z}}$  of real-valued random variables in ${\mathbb L}^p $ and adapted to a non-decreasing filtration $(\F_k)_{k \in {\mathbb Z}}$,  
and let us prove it for $n$. 
Set $a_p =  2 c_p +1$. By the triangle inequality
\begin{equation}\label{extra2}
S_n^* \leq \max_{1\leq k\leq n}
 \left | \sum_{i=1}^k [X_i- \E (X_i| \mathcal F_{i-1})]
 \right | 
 + \max_{1\le k\le n}
 \left | \sum_{i=1}^k \E(X_i | \mathcal F_{i-1})\right |  \, .
\end{equation}
By  von Bahr-Esseen's inequality together with Doob's maximal inequality for martingales,
\begin{equation}\label{ineq1PUWNS}
\left \| \max_{1\le k\le n}
 \left | \sum_{i=1}^k (X_i-\mathbb{E}(X_i|\mathcal F_{i-1}))
 \right | \right \Vert_p
 \leq c_p \left (  \sum_{i=1}^n \Vert X_i-\mathbb{E}(X_i|\mathcal F_{i-1}) \Vert_p^p \right )^{1/p} \leq  2 c_p \left (  \sum_{i=1}^n \Vert X_i \Vert_p^p   \right )^{1/p} \, .  
\end{equation}
To estimate the impact of the second term in the right-hand side of (\ref{extra2}),  we
start by writing $n=2m,$ or $n=2m+1$ according to a value odd or
even of $n$. Notice that
\begin{equation} \label{ineq2PUW}
\left \Vert \max_{1\le k\le n}
 \left |  \sum_{j=1}^k
 \mathbb{E}(X_{i}|\mathcal{F}_{i-1}) \right | \right \Vert_p
 \le \left \Vert \max_{1\leq k\leq m}\left | \sum_{i=1}^{2k}\mathbb{E}
(X_{i}|\mathcal{F}_{i-1})\right | \right \Vert_p 
 +\left \Vert \max_{0\leq k\leq m}
 \left | \mathbb{E} (X_{2k+1}|\mathcal{F}_{2k})\right | 
 \right \Vert_p \, .
\end{equation}
The second term in the right hand side of (\ref{ineq2PUW}) is
estimated in a trivial way:
\begin{equation} \label{ineq3PUWNS}
\left \Vert \max_{0\leq k\leq m}
 |\mathbb{E}(X_{2k+1}|\mathcal{F}_{2k})|\right \Vert_p \leq \left ( \sum_{k=0}^m \Vert \mathbb{E}(X_{2k+1}|\mathcal{F}_{2k}) \Vert_p^p \right )^{1/p}
  \leq \left ( \sum_{i=1}^n \Vert X_i \Vert_p^p \right )^{1/p} \, ,  
\end{equation}
since $m$ is such that $n=2m$ or $n=2m+1$. For the first term in the right hand side of (\ref{ineq2PUW}), we set
\[
Y_i = \E (X_{2i-1} | {\mathcal F}_{2i-2})  + \E (X_{2i} | {\mathcal F}_{2i-1} ) \, , \,  W_i = \sum_{j=1}^i Y_j   \text{ and  } {\mathcal G}_{ i } = {\mathcal F}_{ 2 i -1}\, ,
\] 
and we note that 
\[
 \left \Vert \max_{1\leq k\leq m}\left | \sum_{i=1}^{2k}\mathbb{E}
(X_{i}|\mathcal{F}_{i-1})\right | \right \Vert_p  =  \left \Vert \max_{1\leq k\leq m}\left | \sum_{i=1}^{k} Y_i \right | \right \Vert_p \, .
\]
In addition, $(Y_k)_{k \in {\mathbb Z}}$ is a sequence of real-valued random variables in ${\mathbb L}^p $ and adapted to the non-decreasing filtration $({\mathcal G}_k)_{k \in {\mathbb Z}}$. 
By the induction hypothesis, noticing that $m < 2^{r-1} \leq n$, 
\[
 \left \Vert \max_{1\leq k\leq m}\left | \sum_{i=1}^{k} Y_i \right | \right \Vert_p  \leq a_p \left (  \sum_{j=1}^{m} \Vert Y_j \Vert_p^p \right )^{1/p}  + 2^{(p-1)/p} a_p 
 \sum_{j=0}^{r-2}  \left ( \sum_{k=1}^{2^{r-1 - j}} \Vert {\mathbb E} ( W_{k 2^j } -W_{(k-1) 2^j }  | {\mathcal G}_{(k-1) 2^j} ) \Vert_p^p \right )^{1/p} \, .
\]
But
\[
\Vert {\mathbb E} ( W_{k 2^j } -W_{(k-1) 2^j }  | {\mathcal G}_{(k-1) 2^j} ) \Vert_p \leq \Vert {\mathbb E} ( S_{k 2^{j+1} } -S_{(k-1) 2^{j+1} }  | {\mathcal F}_{(k-1) 2^{j+1}} ) \Vert_p  \, . 
\]
On another hand, 
\[
\sum_{j=1}^{m} \Vert Y_j \Vert_p^p  \leq 2^{p-1}  \sum_{i=1}^{m}  \left ( \Vert \E (X_{2i-1} | {\mathcal F}_{2i-2}) \Vert_p^p + \Vert \E (X_{2i} | {\mathcal F}_{2i-1}) \Vert_p^p \right )  
 \leq  2^{p-1}  \sum_{i=1}^{n}  \Vert \E (X_{i} | {\mathcal F}_{i-1}) \Vert_p^2 \, .
\]
Therefore
\begin{multline*}
 \left \Vert \max_{1\leq k\leq m}\left | \sum_{i=1}^{2k}\mathbb{E}
(X_{i}|\mathcal{F}_{i-1})\right | \right \Vert_p  \leq 2^{(p-1)/p}  a_p \left (   \sum_{i=1}^{n}  \Vert \E (X_{i} | {\mathcal F}_{i-1}) \Vert_p^p   \right )^{1/p} 
\\ +  2^{(p-1)/p} a_p 
 \sum_{j=1}^{r-1}  \left ( \sum_{k=1}^{2^{r-j }} \Vert {\mathbb E} ( S_{k 2^j } -S_{(k-1) 2^j }  | {\mathcal F}_{(k-1) 2^j} ) \Vert_p^2 \right )^{1/2}  \, ,
\end{multline*}
which gives since $n < 2^r$, 
\begin{equation} \label{ineq4PUWNS}
 \left \Vert \max_{1\leq k\leq m}\left | \sum_{i=1}^{2k}\mathbb{E}
(X_{i}|\mathcal{F}_{i-1})\right | \right \Vert_p   \leq  2^{(p-1)/p} a_p 
 \sum_{j=0}^{r-1}  \left ( \sum_{k=1}^{2^{r-j }} \Vert {\mathbb E} ( S_{k 2^j } -S_{(k-1) 2^j }  | {\mathcal F}_{(k-1) 2^j} ) \Vert_p^2 \right )^{1/2}  \, .
\end{equation}
So, overall, starting from \eqref{extra2} and taking into account the upper bounds \eqref{ineq1PUWNS}, \eqref{ineq2PUW}, \eqref{ineq3PUWNS} and \eqref{ineq4PUWNS}, Inequality 
\eqref{eq:bs0314NSpmoinsque2} follows proving the induction hypothesis at step $n$. $\Diamond$

\subsection{Proof of Proposition \ref{weakvbe}}
 Let $M>0$. For every $1\le k\le n$, define 
$$
\tilde D_k:= D_k{\bf 1}_{\{|D_k|\le y \}} - 
\E \left (D_k{\bf 1}_{\{|D_k|\le y \}}|\F_{k-1} \right )\, ,
$$ 
so that $(\tilde D_k)_{1\le k\le n}$ is a sequence of martingale diferences. 
We have 

\begin{equation*}
\P\left (\max_{1\le k\le n}\left |\sum_{j=1}^k (D_j-\tilde D_j) \right |\ge y/2 \right ) 
\le \frac{4 }{y} \sum_{k=1}^n\E \left (|D_k|{\bf 1}_{\{|D_k|>y \}} \right )\, .
\end{equation*}
Now, 
\begin{equation*}
\E\left (|D_k|{\bf 1}_{\{|D_k|>y)\}} \right )
=\int_0^y \P(|D_k|>y) dt +\int_y^{+\infty} \P(|D_k|>t) dt
\le \frac{p}{p-1} \|D_k\|_{p,\infty}^py^{1-p} \, .
\end{equation*}
Hence, 
\begin{equation}\label{B1}
\P \left (\max_{1\le k\le n} \left |\sum_{j=1}^k (D_j-\tilde D_j) \right |\ge y/2 \right ) 
\le   \frac{4p}{y^p(p-1)}\sum_{k=1}^n \|D_k\|_{p,\infty}^p\, .
\end{equation}

On another hand, by Doob's maximal inequality, 
\begin{equation*}
\P \left (\max_{1\le k\le n} \left |\sum_{j=1}^k \tilde D_j \right |\ge y/2 \right ) 
\le  \frac{4}{y^2}  \sum_{k=1}^n \E(D_k^2{\bf 1}_{\{|D_k|\le y\}})\, .
\end{equation*}
Now,
\begin{equation*}
\E \left (D_k^2{\bf 1}_{\{|D_k|\le y\}} \right )\le
 \int_0^y 2t \P(|D_k|>t) dt \le \frac{2}{2-p}\|D_j\|_{p,\infty}^p y^{2-p}\, .
\end{equation*}
Hence,
\begin{equation}\label{B2}
\P \left (\max_{1\le k\le n} \left |\sum_{j=1}^k \tilde D_j \right |\ge y/2 \right ) 
\le   \frac{8}{y^p(2-p)} \sum_{k=1}^n \|D_k\|_{p,\infty}^p \, . 
\end{equation}
The result follows from \eqref{B1} and \eqref{B2}.

\section{General cocycles}\label{Sec:GC}

\setcounter{equation}{0}

It turns out that all the results obtained under moments greater than 2 
made use of a martingale-coboundary decomposition with bounded 
(in ${\mathbb L}^\infty$) coboundary and of the fact that we study partial sums associated with a cocycle. Another ingredient of general nature used in the proofs is Proposition 3.1 of Benoist-Quint \cite{BQ}. In particular all the results   obtained under moments greater than 2 may be generalized to cocycles admitting such a martingale-coboundary decomposition. Such cocycles are called \emph{centerable} in \cite{BQ}. 

\smallskip

We shall also give sufficient conditions under which the results under 
moments weaker than 2 holds for general cocycles. 

\medskip

Let us describe the situations that should be considered in the sequel.

\medskip

Let $G$ be a locally compact second countable group. Let $X$ be  compact 
and second countable. Assume that $G$ acts continuously on $X$ and denote 
that action by $g\cdot x$. 

\smallskip

Let $\sigma\,:\, G\times X\to \R$ be a cocycle, meaning that it satisfies  the equality \eqref{cocycle-prop}.
We shall only be concerned with continuous cocycles. Given a continuous 
cocycle, define $\sigmas(g):=\sup_{u\in X}|\sigma(g,u)|$ for every 
$g\in G$.

\medskip

Let $\mu$ be a probability measure on the Borel sets of $G$.  

\medskip

Assume that there exists a \emph{unique} $\mu$-invariant probability $\nu$
on the Borel sets of $X$, that is a unique probability satisfying 
\eqref{inv22-07}.

\medskip

Let $(\Omega,\F,\P)$ be a probability space. Assume that there exists a sequence
$(Y_n)_{n\ge 1}$of  iid random variables on $(\Omega,\F,\P)$ taking values 
in $G$ with common law $\mu$. Define $A_n:=Y_n\cdots Y_1$ for every      
$n\ge 1$ and $A_0={\bf e}$ the neutral element of $G$. 

\medskip

Our goal is to study the sequence defined by 
$$
S_{n,u}:= \sigma(A_n,u)=\sum_{k=0}^{n-1}\sigma(Y_{k+1}, A_k\cdot u)
\qquad \forall n\ge 1,\, \forall u\in X\, .
$$

\begin{Definition}
We say that $\sigma$ is \emph{centerable} if $\sigmas\in  {\mathbb L}^1(\mu)$ and if there exist a cocycle $\sigma_0$
 and a continuous function $\psi$ on $X$ such that 
$\int_G\sigma_0(g,u)\mu(dg)=\lambda_\mu$ for every $u\in X$, where $\lambda_\mu:=\int_{G\times X}\sigma(h,v)\mu(dh)\nu(dv)$, and 
\begin{equation}\label{centerable}
\sigma(g,u)
=\sigma_0(g,u)+\psi(u)-\psi(g\cdot u)\qquad \forall (g,u)\in G\times X\, .
\end{equation}

\end{Definition}
\begin{Remark}
A sufficient condition for $\sigma$ to be centerable is Gordin's condition: 
$$
\sum_{n\ge 0} \sup_{u\in X} |\E(\sigma(Y_{n+1}, A_n\cdot u))-\lambda_\mu|
= \sum_{n\ge 0} \sup_{u\in X} \left|\int_{G\times G}\sigma(g,g'\cdot u)
\mu(dg)\mu^{*n}(dg)\, -\lambda_\mu\right|<\infty\, .
$$
\end{Remark}

\medskip

$\bullet$ Assume that there exist $r>0$ and $\delta >0$ such that 
$$
\int_G {\rm e}^{\delta\,  \sigmas^r(g)}\mu(dg)<\infty\, .
$$
If $\sigma$ is centerable, then the conclusions of Theorem 
\ref{LDsuperexp} and Theorem 
\ref{LDsousexp} hold with $S_{k,u}$ in place of $\log\|A_k x\|$ 
for the corresponding value of $r>0$.

\medskip

$\bullet$ Assume $\sigmas\in  {\mathbb L}^2(\mu)$ and 
 \begin{equation}
 \limsup_{n \rightarrow \infty} 
  \frac{n}{b_n^2}\log n \mu \left \{ \sigmas > b_n \right \} 
  = - \infty \, , 
\end{equation}
for some sequence $(b_n)_{n\ge 1}$ satisfying \eqref{condseq}. 
Then, if $\sigma$ is centerable, the conclusion of Theorem \ref{MDPth} 
holds with $S_{[nt],u}$ in place of $\log\|A_{[nt]}x\|$.

\medskip

In the same way, if $\sigma$ is centerable, in the case of weak moments of order $p >2$ (resp.  strong moments of order $p \geq 2$), the conclusion of  Theorem \ref{devweak}  (resp. of Theorem \ref{devstrong}) holds 
with $S_{k,u}$ in place of $\log\|A_k x\|$. 

\medskip

Let us now mention results under weak moments of order $p$, $1<p<2$ or under 
(strong) moments of order $1\le p<2$.

\medskip

$\bullet$ Let $1<p<2$. Assume that 
$$
\sup_{t>0} t^p\mu\{\sigmas >t\}<\infty\, ,
$$
and that
$$
\sum_{n\ge 1} n^{-1-1/p}\sup_{u\in X} \left|\E(\sigma(Y_n,A_{n-1}\cdot u))\, -
\lambda_\mu\right|<\infty\, .
$$
Then the conclusion of Theorem \ref{devweak} holds with $S_{k,u}$ 
in place of $\log\|A_k x\|$.

\medskip

$\bullet$ Let $1\le p<2$. Assume that $\sigmas\in {\mathbb L}^p(\mu)$
and that there exists $q>p$ such that
$$
\sum_{n\ge 1} n^{-1-1/q}\sup_{u\in X} \left|\E(\sigma(Y_n,A_{n-1}\cdot u))\, -
\lambda_\mu\right|<\infty\, .
$$
Then the conclusion of Theorem \ref{devstrong} holds with $S_{k,u}$ 
in place of $\log\|A_k x\|$.

\end{document}